\documentclass[11pt]{elsarticle} 
 
\usepackage{amsmath}
\usepackage{amssymb}
\usepackage[all]{xy} 
\usepackage{enumerate} 
\usepackage{array}
\usepackage{amsthm}
\usepackage{amscd}
\usepackage{color}  
\usepackage{hyperref} 
\hypersetup{
    colorlinks, 
    citecolor=red,
    filecolor=black, 
    linkcolor=blue,
    urlcolor=black
}
\hypersetup{linktocpage} 
\usepackage{multirow}

\def \D{\mathbb{D} \,}
\def \R{\mathbb{R} \,}

\newcommand \N{\mathbb{N}}

\newdefinition{rmk}{Remark}

\journal{Applied Mathematics and Computation}

\begin{document}
\begin{frontmatter}

\title{Optimal sampling patterns for Zernike polynomials}


\author[address1]{D. Ramos-L\'opez}
\ead{drl012@ual.es}
\author[address1]{Miguel \'Angel S\'anchez-Granero}
\ead{misanche@ual.es}
\author[address3]{Manuel Fern\'andez-Mart\'inez}
\ead{fmm124@gmail.com}
\author[address1,address2]{A. Mart\'{\i}nez--Finkelshtein\corref{cor1}}
\ead{andrei@ual.es}
\cortext[cor1]{Corresponding author.}
\address[address1]{Department of Mathematics,
University of Almer\'{\i}a, Spain}
\address[address2]{Instituto Carlos I de F\'{\i}sica Te\'{o}rica y Computacional,
Granada University, Spain}
\address[address3]{University Center of Defense at the Spanish Air Force Academy, MDE-UPCT, Santiago de la Ribera, Murcia, Spain}

\begin{abstract}
A pattern of interpolation nodes on the disk is studied, for which the interpolation problem is theoretically unisolvent, and which renders a minimal numerical condition for the collocation matrix when the standard basis of Zernike polynomials is used. It is shown that these nodes have an excellent performance also from several alternative points of view,  providing a numerically stable surface reconstruction, starting from both the elevation and the slope data. Sampling at these nodes allows for 
a more precise recovery of the coefficients in the Zernike expansion of a wavefront or of an optical surface. 
\end{abstract}

\begin{keyword}
Interpolation \sep Numerical condition \sep Zernike polynomials  \sep Lebesgue constants
\end{keyword}

\end{frontmatter}

\section{Introduction} \label{sec1}

 Zernike or circular polynomials \cite{ZER34} constitute a set of basis functions, very popular in optics and in optical engineering,  especially appropriate to express wavefront data due to their connection with classical aberrations. Some of their applications in optics include optical engineering \cite{BOR70}, aberrometry of the human eye \cite{CAR05,CON83,KLY04}, corneal surface modeling \cite{ISK01,SMO03,SMO05} and other topics \cite{ARE06,BRA02,JAN02,NAM09} in optics and ophthalmology. Due to their optical properties and pervasiveness, Zernike polynomials are included in the ANSI standard to report eye aberrations \cite{ANSI04}.

In almost every practical application, an optical surface or wavefront is sampled at a finite set of points, followed by a fit of the collected data by a linear combination of the Zernike basis with the purpose to determine the corresponding coefficients (sometimes called modes). This is normally done by means of the standard technique of linear least squares, which reduces to interpolation when the size of the data and the dimension of the basis match. 
The process is often ill-conditioned, as its stability strongly depends on an adequate choice of the sampling nodes. Analogous situation arises, for instance, when the Zernike basis is used to fit the slopes of the wavefront in a Shack-Hartmann device \cite{Lane}. Their partial derivatives are used in a least squares fit and if the sampled nodes are not chosen carefully, the resulting Zernike modes might be totally inaccurate.

Different sampling patterns can be found in \cite{Diaz-Santana2005,Navarro2009,Navarro2011}. Some of them are based on  random or pseudo-random points drawn according to a probability distribution on the disk. Other schemes include regular or quasi-regular grids, such as squared or hexagonal Cartesian grids, regular polar grids or hexapolar grids, that cover  the surface of the disk more or less uniformly. However, these sampling patterns generally produce an ill-conditioned collocation matrix even for moderate Zernike polynomial orders. This issue was addressed in \cite{Navarro2009} for the reconstruction of the wavefront from its slopes, putting forward a spiral arrangement as the best best-performing sampling pattern for this problem. However, the numerical results show that even this pattern is not totally satisfactory for the elevation data.  

The main goal of this paper is to discuss in a certain sense optimal patterns for sampling and interpolation on the disk using the basis of the Zernike polynomials for moderate degrees, used in practical applications in optics and ophthalmology. They render well-conditioned collocation matrices and provide numerically stable surface reconstruction, starting from both the elevation and the slope data.


\section{Methods} \label{sec2}

\subsection*{Zernike polynomials}

Zernike polynomials are usually defined in polar coordinates using the double-index notation (see e.g.~\cite[Ch.~IX]{BOR70}),
\begin{equation}
Z_n^m(\rho,\theta) = 
\begin{cases} 
\gamma_n^m R_n^{|m|}(\rho) \cos(m\theta), &\text{if } m\geq 0, \\ 
\gamma_n^m R_n^{|m|}(\rho) \sin(|m|\theta), &\text{if } m< 0,
\end{cases}  
\end{equation}
where $n\geq 0$, $|m|\leq  n $, and $n-m$ is  even, $\gamma_n^m$ are normalization constants, and the radial part $R_n^m$  is 
$$ 
R_n^{|m|}(\rho)=\sum_{s=0}^{(n-|m|)/2}\frac{(-1)^s (n-s)!}{s!((n+|m|)/2-s)! ((n-|m|)/2-s)!}\rho^{n-2s}, 
$$
which can be expressed in terms of shifted Jacobi polynomials $P_n^{(0, m-n)}$. Alternatively, a single index notation is used, and the conversion from $Z_n^m$ to $Z_j$ is made by the formula \cite{ANSI04}
\begin{equation}
\label{formulaJ}
j =\frac{n (n+ 2) + m}{2}\in \N \cup \{0\}.
\end{equation}

Functions $Z_n^m$ are actually polynomials in Cartesian coordinates $(x,y)$. Index $n$ is called the \emph{radial order} of $Z_n^m$. It is easy to see that the number of distinct Zernike polynomials of radial order $\leq n$ is
\begin{equation}
\label{defN}
N=\frac{(n+1)(n+2)}{2},
\end{equation}
which matches the dimension of algebraic polynomials in two variables of total degree $\leq n$. In fact, Zernike polynomials are a complete polynomial set on the unit disk $\D=\{(x,y)\in \R^2: x^2+y^2\leq 1\}$,  orthogonal with respect to the area measure on the disk. 

The normalization factor $\gamma_n^m$ for simplicity can be set to $1$;  however,  with $\gamma_n^m = \sqrt{ (2-\delta_{0,m})(n+1)}$, where $\delta$ is the Kronecker delta, the set of polynomials is orthonormal:
$$
\iint_{\D}Z_n^m(\rho,\theta)Z_r^s(\rho,\theta)\rho \,d\rho d\theta=\delta_{n,r}\delta_{m,s}.
$$
In what follows, notation $Z_n^m$ stands for the orthonormal Zernike polynomials. 

The standard Fourier theory can be easily extended to the circular polynomials. In particular, any function $W\in L^2$ defined over $\D$ can be represented as
\begin{equation}
\label{Fourierseries}
W(\rho,\theta) = \sum_{m,n} c_{n}^m\, Z_n^m(\rho,\theta), \quad c_{n}^m = \iint_{\D} W(\rho,\theta)  Z_n^m(\rho,\theta) \rho\,  d\rho d\theta.
\end{equation}
However, in real applications function $W$ is sampled only at a discrete (and finite) set of nodes, and its Zernike (Fourier) coefficients $ c_{n}^m$ must be recovered using only this available information. A common and statistically meaningful procedure seeks a solution in the (weighted) least squares sense to a linear system whose matrix (known as the \emph{collocation matrix}) consists of evaluations of the Zernike basis in the given set of points or nodes (see e.g.~\cite{ISK02}). In the case when the number of nodes matches the dimension of the polynomial subspace used for approximation (called ``critical sampling'' in \cite{Navarro2009}), the problem boils down to the polynomial interpolation of the function $W$ at a given set of nodes.

 It is worth mentioning an alternative approach, called \emph{hyperinterpolation} \cite{hyperinterpolation}, where $W$ is approximated 
 by truncations of the series in \eqref{Fourierseries}, but the Fourier coefficients $c_{n}^m$ are computed using quadrature formulas evaluated at the discrete set of nodes. The analysis of this method lies beyond the scope of this paper. 

\subsection*{Goodness of the sampling patterns}

It is well known that 
 in the multivariate case the unisolvence of the interpolation problem for arbitrary  nodes is not guaranteed, so not every sampling pattern is acceptable (several configurations of interpolation points on the disk that guarantee unisolvence are well-known and can be found in the literature, see e.g.~\cite{Sauer96,Bojanov2003,Xu2004,Cuyt2012}). Moreover, the error in approximating a function by its interpolating polynomial depends on the interpolation nodes: standard upper bounds for the error are based on the so-called \emph{Lebesgue constants} corresponding to these nodes, which give the norm of the interpolation as a projection operator onto the polynomial subspace (see e.g.~\cite{Cuyt2012,Phillips2003}). 
 
Thus, if we are interested in the set of interpolation points which gives the smallest possible upper bound
on the interpolation error in an arbitrary continuous function, an optimal choice 
of interpolation points (at least, in this sense) is given by those which minimize the Lebesgue constant. 
The asymptotic theory of these interpolation sets is rather well understood. For instance,   as it was established in \cite{Sunder84}, the order of growth of the Lebesgue constants on the disk for algebraic polynomials  of total degree $\leq n$ is $\geq \mathcal O(\sqrt{n})$. 
Moreover, the sub-exponential growth of the Lebesgue constants implies, among other facts, the weak-* convergence of the nodes counting measure  to the (pluripotential theory) equilibrium measure of the disk,  given by 
$(2\pi \sqrt{1-x^2-y^2})^{-1}dxdy$ (see \cite{bloom92} and \cite{bloom12}).

Closely related is the notion of \emph{admissible} and  \emph{weakly admissible meshes}  (see \cite{Calvi08} and \cite{Piazzon}), namely sequence of discrete subsets $\mathcal A_n$ of a compact set $K$ such that
$$
\sup_p \, \frac{\max \{| p(x)|: x\in K\}}{\max \{| p(x)|: x\in \mathcal A_n\}}
$$
over the polynomials of degree $\leq n$ is either uniformly bounded (admissible) or grows at most polynomially in $n$ (weakly admissible).

Minimizing the Lebesgue constant amounts to solving a large scale non-linear optimization problem, to which the true solution is not explicitly known, even in the case of univariate interpolation. In fact, no explicit examples of multivariate interpolation sets providing an at most polynomial rate of growth of the Lebesgue constants are currently available. In the case of the square $[0,1]^2$, the most popular nodes are given by the so-called Padua points \cite{Padua}, but their construction seems to be hardly generalizable to other sets.  For the disk, no explicit configurations of interpolation points obeying the order of growth of $  \mathcal O(\sqrt{n})$ are known. Good candidates are points minimizing certain energy on the disk, such as Leja (giving a polynomial growth of the Lebesgue constants, as established in \cite{Leja}) or Fekete points \cite{Briani}.  The so-called \emph{Bos arrays} (see \cite{bos81}, as well as \cite{bloom92} and \cite[Section 8]{bloom12}) provide a polynomial growth of the Lebesgue constants, which usually suffices for practical applications. 
The Lebesgue constants  for several other unisolvent configurations have been numerically analyzed in \cite{Cuyt2012}.

However, these criteria of optimality of nodes do not take into account the numerical aspects of the interpolation. For instance, if we have a polynomial basis fixed, interpolation boils down to finding the coefficients of the expansion of the interpolating polynomial in terms of this basis, which is equivalent to solving a linear system given by the so-called collocation matrix. It is well know that  the actual choice of the basis greatly influences the accuracy of the solution to this problem. 

This is the situation with the actual practice in ophthalmology and the visual science: the polynomial basis is usually given a priori by the orthonormal Zernike polynomials (see e.g.~\cite{ISK02} for other bases). Assuming this as the starting hypothesis, we analyze the dependence of the numerical stability of the associated linear system from the selection of the nodes. 

 In theory, unisolvence of the interpolation problem is equivalent to the regularity of the corresponding collocation matrix.  
However, from the practical point of view invertibility of the collocation matrix is not sufficient, since ill-conditioned problems are numerically infeasible. The numerical conditioning of a system of linear equations can be measured by the \emph{condition number} $\kappa(A)$ of the system matrix $A$ with respect to inversion, see e.g.~\cite{Trefethen97}. Roughly speaking, $\kappa(A)\approx 10^s$ means a possible loss of about $s$ digits of accuracy in the solution of the linear system. In particular, when working in the IEEE double precision  with matrices with condition numbers of order  $10^{16}$, all the significant digits may be lost. 

Our goal is to put forward a construction of a Bos array of interpolation nodes on the unit disk for which the corresponding collocation matrix built from Zernike polynomials is particularly well-conditioned, at least for moderate degrees. Motivated essentially by basic applications in visual sciences, we do not tackle the asymptotic problem here.

\begin{figure}
\centering 
\begin{tabular}{cc}
\includegraphics[width=0.45 \textwidth]{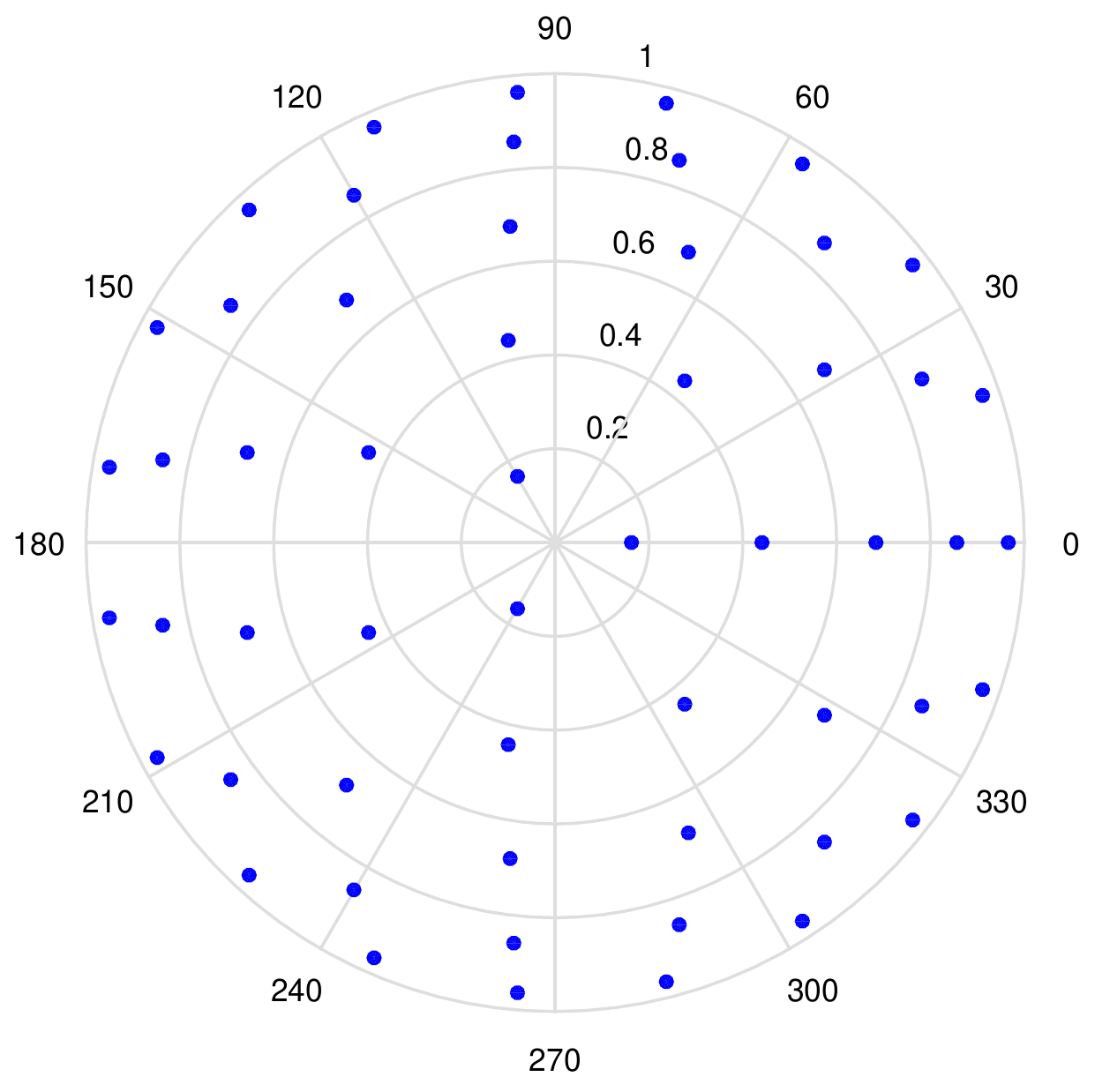} &  
\includegraphics[width=0.45 \textwidth]{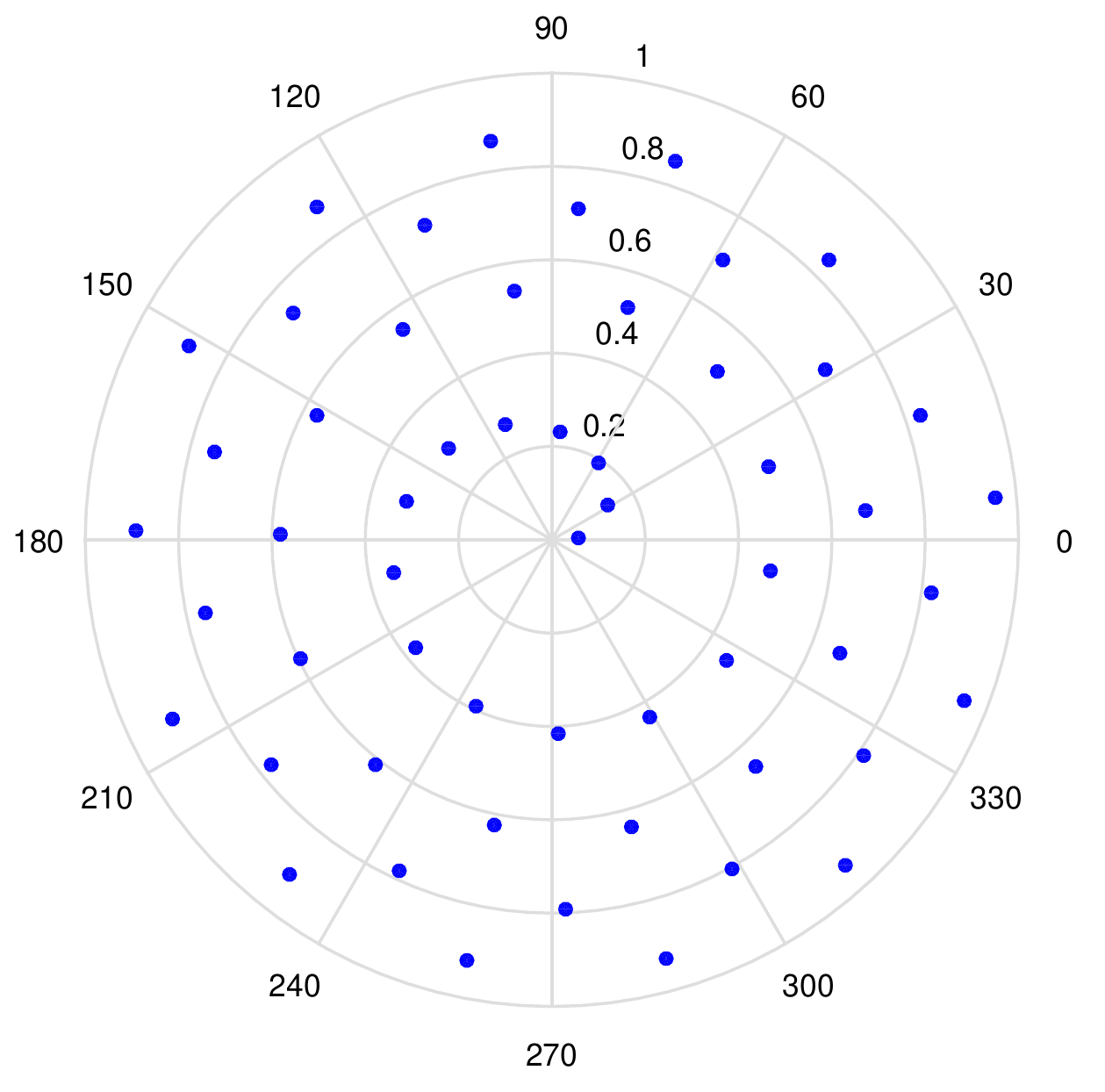} 
\\ \includegraphics[width=0.45 \textwidth]{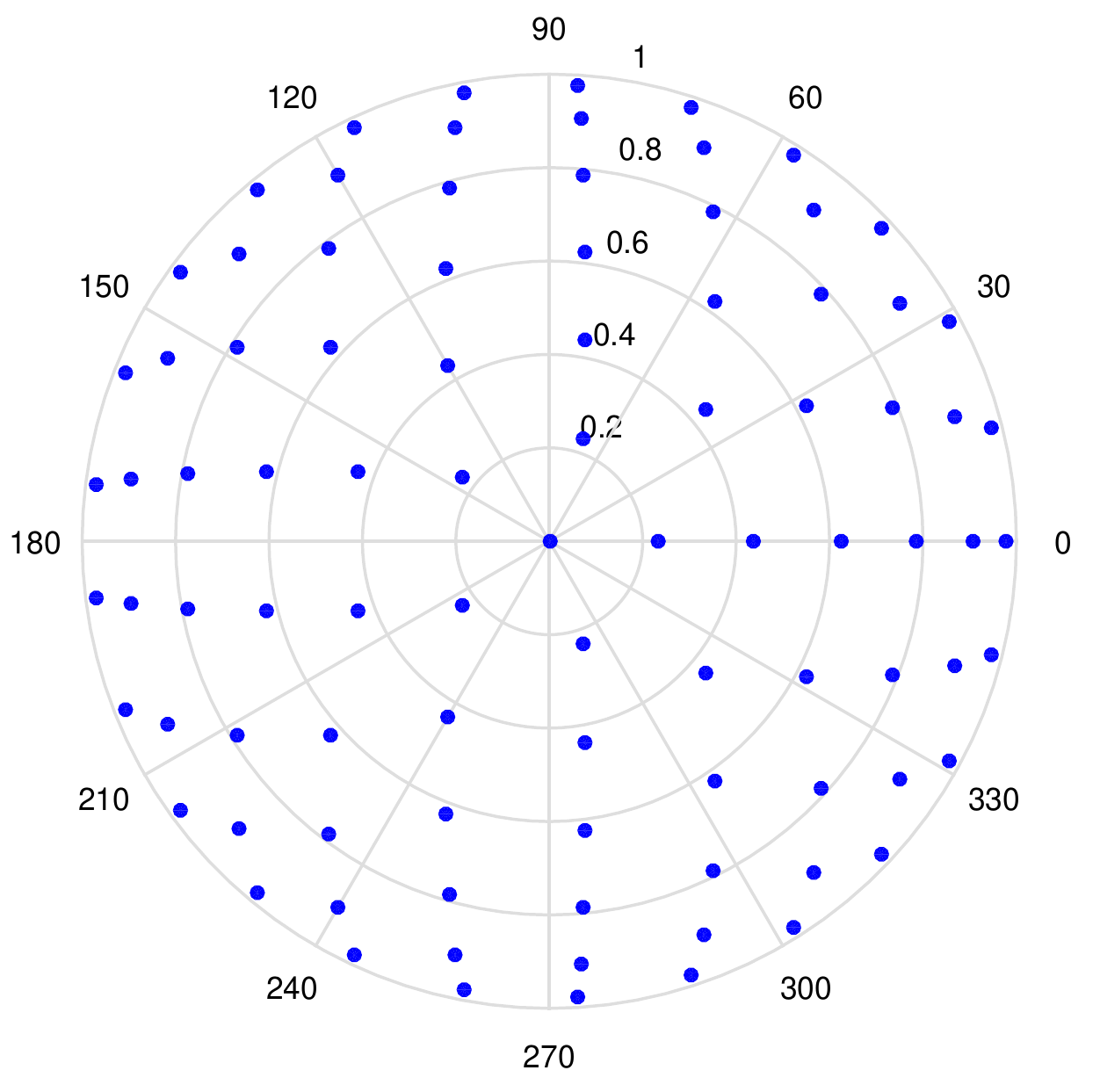} &  
\includegraphics[width=0.45 \textwidth]{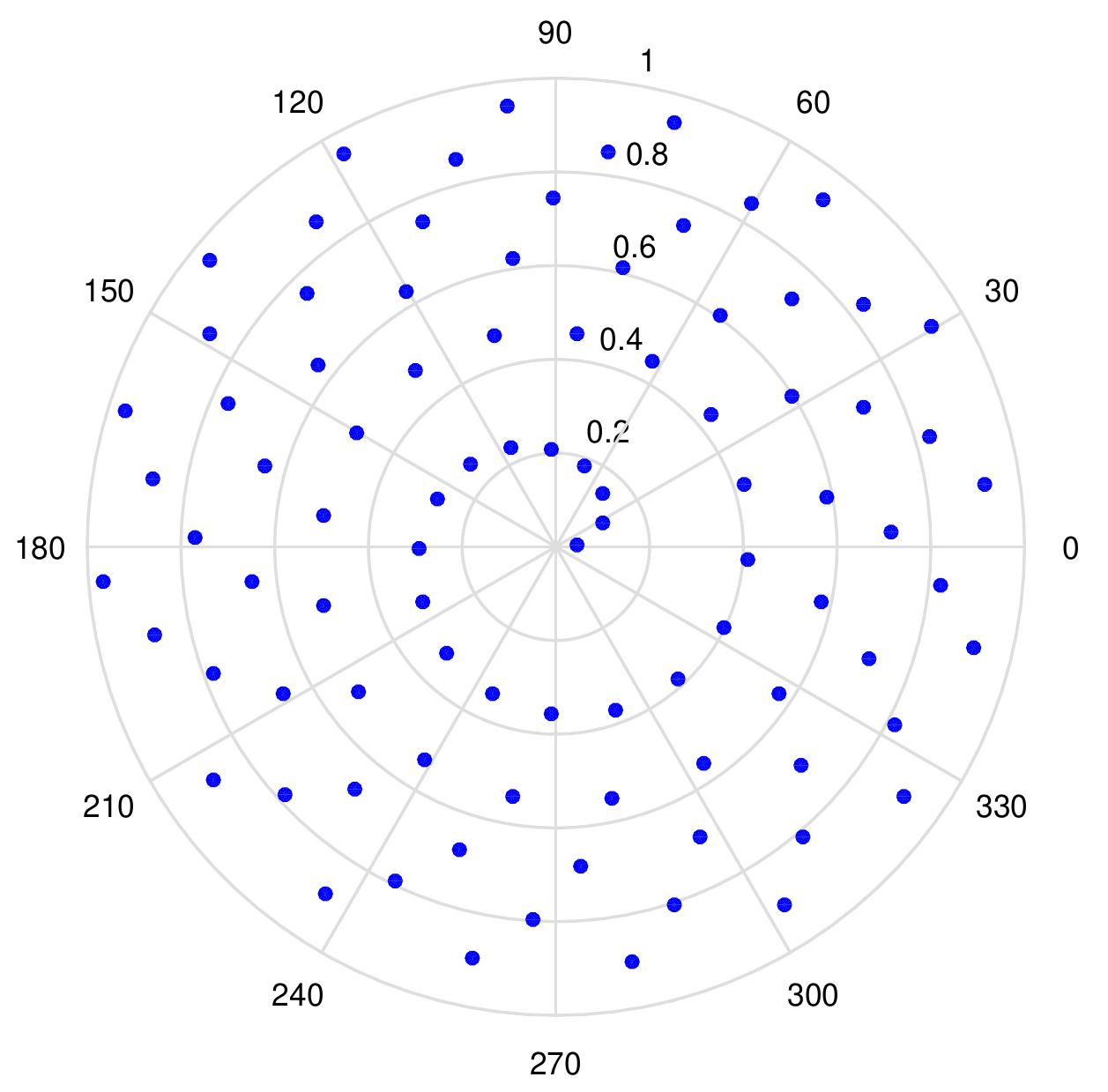}
\end{tabular}
 \caption {Concentric sampling patterns for Zernike polynomials (left) and spiral sampling pattern, as defined in \cite{Navarro2009} (right), for Zernike polynomials of radial order 9 (top) and 12 (bottom).}
\label{figure:newPatterns}
\end{figure}

We will also make a  comparison with the recently introduced spiral sampling \cite{Navarro2009} for the reconstruction of the wavefront from its slopes, which allowes for a stable use of the Zernike polynomials up to radial order 15, approximately (see Figure~\ref{figure:newPatterns}, right, where the spiral pattern is  depicted for radial orders 9 and 12, respectively). 

\subsection*{Interpolation nodes with low condition numbers}

We  minimize the condition number of the collocation matrix for a pattern of nodes on concentric  circles (the \emph{Bos arrays}), for which unisolvence has been established (see \cite{bos81}, as well as \cite{Bojanov2003}). Namely, given the maximal radial order $n$ of Zernike polynomials, we choose 
$$
1 \geq r_1> r_2 >\dots > r_k\geq 0,  \quad k=k(n) = \left\lfloor \frac{n}{2}\right \rfloor+1
$$ 
radii (where  $\lfloor \cdot \rfloor$ is the floor operator); on the $i$-th circle with center at the origin and radius $r_i$ we place
$$
n_i = 2 n+5-4i
$$
equally spaced nodes. Notice that
$$
\sum_{i=1}^k n_i = \frac{(n+1)(n+2)}{2} = N,
$$
namely, the number of Zernike polynomials of radial order $\leq n$. With this rule, the $k$-th (innermost) circle contains $1$ node if $n$ is even, and $3$ nodes otherwise (see Figure~\ref{figure:newPatterns}, left). Observe that we have not prescribed the exact position of the nodes (only that they are equally spaced). The latter condition is relevant: there are examples in \cite{Bojanov2002} showing that the interpolation problem is not poised for arbitrary points on the circle. 

For the sake of precision, let us set the nodes on the circle of radius $r_i$ having the arguments $2\pi t/n_i$, $t=0, \dots, n_i-1$. However, we will see that the actual positions of the nodes on the respective circles have a relative relevance. 

If we denote by $P_i$, $i=1, \dots, N$, the nodes constructed according to this rule, and if $Z_j$ are the Zernike polynomials enumerated using formula \eqref{formulaJ}, then the collocation matrix takes the form
\begin{equation}
\label{defA_n}
A_n=\left( Z_{j-1}(P_i)\right)_{i,j=1}^N,
\end{equation}
and we are interested in its condition number $\kappa_2(A_n)$, which can be computed as the ratio of the largest to the smallest singular values of $A_n$ (see \cite{Trefethen97}). Since the singular values are invariant by permutation of rows and columns of $A_n$,  the value of $\kappa_2(A_n)$ is independent of the order of the nodes $P_i$ and of the Zernike polynomials $Z_j$.

In the proposed scheme, the relevant parameters at hand are the radii $r_i$. We consider the problem of
\begin{equation}
\label{minimization}
\min \{\kappa_2(A_n):\, 1 > r_1> r_2 >\dots > r_k\geq 0\};
\end{equation}
the optimal values of $r_j$'s were found using simulated annealing followed by the algorithms for non-linear optimization implemented in the Optimization Toolbox of Matlab, achieving at least 4 digits of precision, which is sufficient for interpolation with Zernike polynomials of radial order $\leq 30$. The constraint $r_1<1$ is for a practical reason: in actual measurements the data in the periphery is usually much less reliable.   

A similar problem, but minimizing the Lebesgue constants for the nodes $\{P_i\}$, was considered in \cite{Cuyt2012,Briani,Carnicer,Gunzburger}. In \cite{Cuyt2012},  the radii were studied in connection with the notion of spherical orthogonality for multivariate polynomials, and values of $r_j$'s, given by the zeros of certain Gegenbauer orthogonal polynomials, were analyzed experimentally.    
In~\cite{Carnicer}, along with the Lebesgue constants, the condition numbers $\kappa_\infty(A_n)=\|A_n\|_\infty \|A_n^{-1}\|_\infty$ were used to optimize the radii.
The nodes described below, found as the solution to the problem \eqref{minimization},  differ from those  obtained in \cite{Cuyt2012} and \cite{Carnicer}. 
 
So far, we have not been able to associate the optimal values of $r_j$'s with any known set. Instead, we will prescribe the quasi-optimal radii $r_j$ 
using the following formula, obtained by the least square fitting of the optimal radii:
\begin{equation}
\label{recCheb}
r_j=r_j(n) = 1.1565 \, \zeta_{j,n} - 0.76535 \, \zeta_{j,n}^2 + 0.60517 \, \zeta_{j,n}^3, \quad 
\end{equation}
where $\zeta_{j,n}$ are zeros of the $(n+1)$-st Chebyshev polynomial of the first kind,
$$
\zeta_{j,n} = \cos\left( \frac{(2j-1) \pi}{2(n+1) }\right), \quad j=1, \dots , k = \left\lfloor \frac{n}{2}\right \rfloor+1.
$$
Given these radii, the interpolation nodes $P_j$ are defined in polar coordinates as
\begin{equation}
\label{nodes}
\left( r_j, 2\pi \frac{s_j-1}{n_j} \right), \quad j=1, \dots, k(n), \quad s_j=1, \dots, n_j.
\end{equation}


We have also studied numerically the behavior of the Lebesgue constants corresponding to our interpolation nodes $\{P_i\}$. Recall \cite{Phillips2003} that the Lebesgue constant $\Lambda_n$ corresponding to the polynomial interpolation with total degree $\leq n$ is the maximum over the disk $\D$ of the function
\begin{equation}
\label{lebesgue1}
\ell(x,y)=\sum_{i=1}^N |\ell_i(x,y)|,
\end{equation}
where $\ell_i(x,y) $ is the $i$-th basic Lagrange polynomial characterized by
$ 
\ell_i(P_j) =\delta_{ij}$. In particular,
\begin{equation}
\label{lebesgue2}
\ell_i(x,y)= \frac{1}{\det A_n}\, \det \left( A_n^{(i)}(x,y)\right),
\end{equation}
where $A_n$ was defined in \eqref{defA_n}, and $A_n^{(i)}(x,y)$ is obtained from $A_n$ by replacing the $i$-th row with $\left( Z_{j-1}(x,y)\right)_{j=1}^N$. 

It is worth mentioning that other nodes can be obtained following \cite{vanbarel,Narayan}, where greedy optimization algorithms to compute a ``good'' set of nodes for multivariate polynomial interpolation are used.

\section{Experimental results} \label{sec3}

We have carried out numerical results in order to assess the performance of several sampling patterns. The random, hexagonal, hexapolar, square and the spiral pattern proposed by Navarro et al. were compared in \cite{Navarro2009}. The results of our analysis corroborate the main conclusion of \cite{Navarro2009}, i.e., that the spiral sampling outperforms the rest of the patterns discussed therein. However, as it could be inferred already from \cite{Carnicer}, the concentric configurations can give better results, at least in terms of the numerical conditioning of the collocation matrix $A_n$. 
In this paper, we compare the spiral sampling pattern from \cite{Navarro2009} with the interpolation scheme, given by \eqref{recCheb}--\eqref{nodes}, and to which we refer as \emph{the optimal concentric sampling}, or the OCS, for short.  

\begin{figure}
\centering \begin{tabular}{lr}
\hspace{-10mm} \includegraphics[width=0.5 \textwidth]{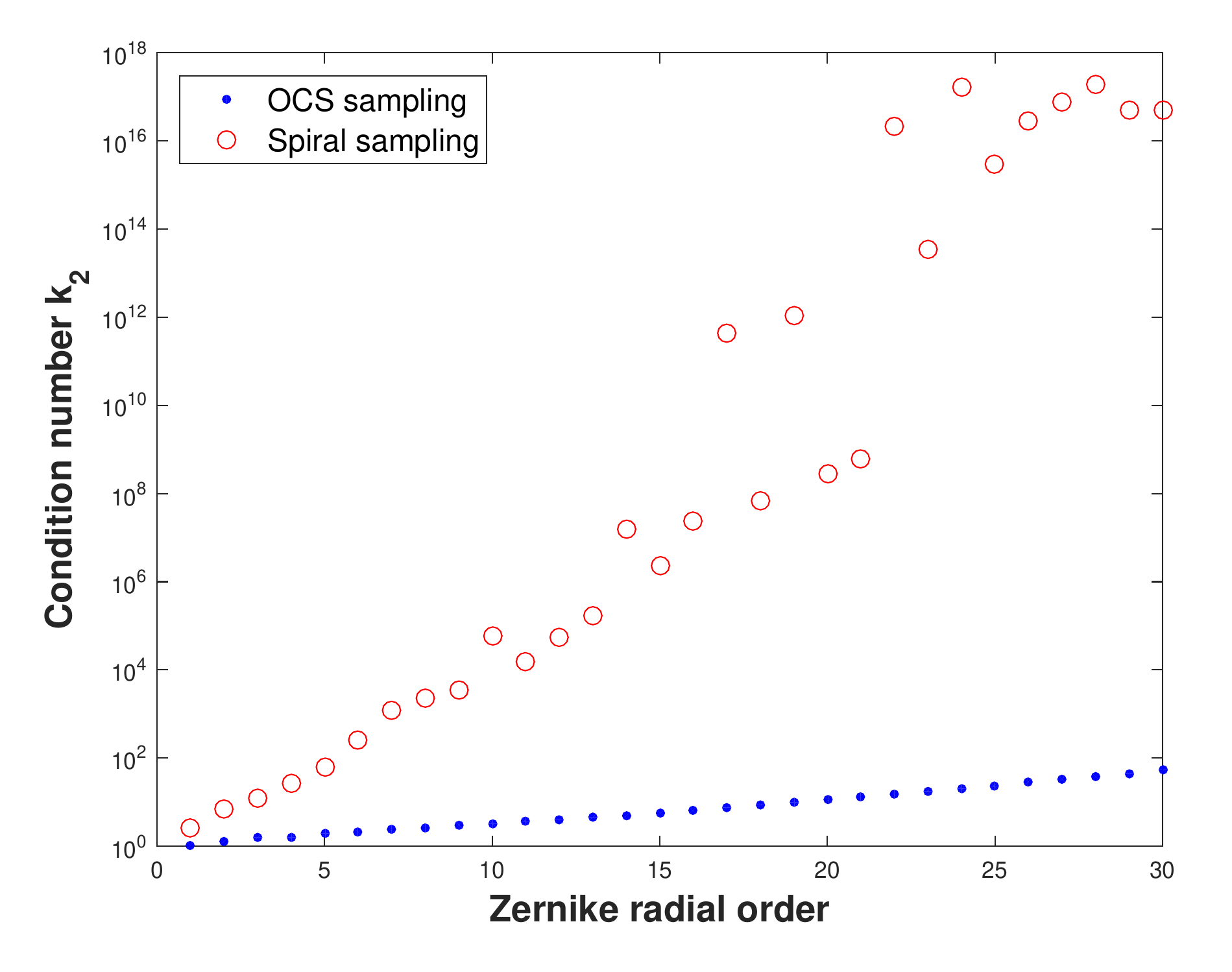} &
\hspace{-5mm} \includegraphics[width=0.5 \textwidth]{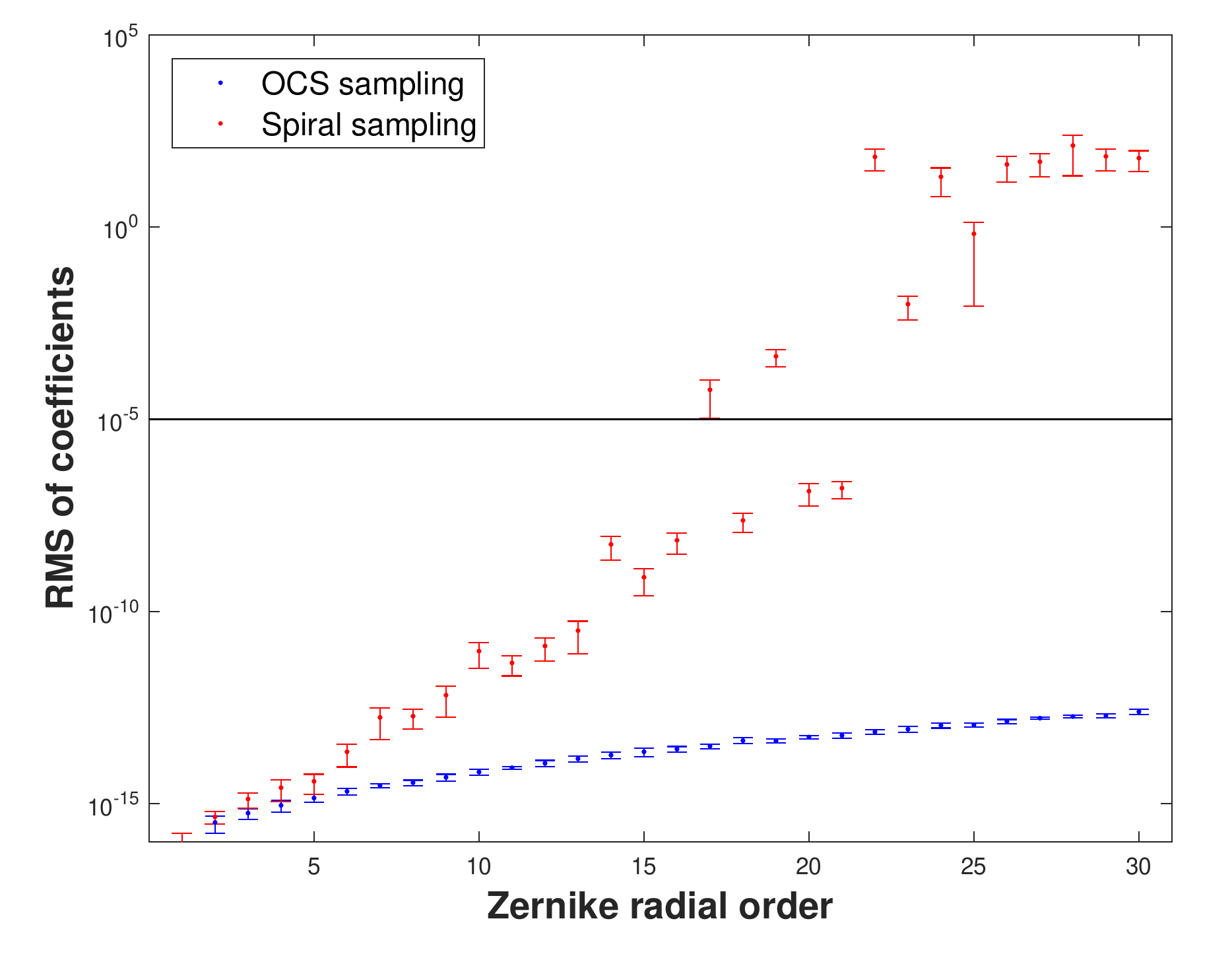}
\end{tabular}
 \caption{Left: condition number $\kappa_2$ of the collocation matrix \eqref{defA_n} as a function of the radial order, for both the spiral pattern and the OCS. Right: range plot of the RMS error of the recovered Zernike modes with respect to the exact coefficients.}
\label{figure:Conditioning}
\end{figure}

To make a proper comparison between the spiral sampling and the OCS, a variety of experiments were carried out. In a wide range of Zernike radial orders, from 1 to 30, the condition numbers $\kappa_2$  of the collocation matrix \eqref{defA_n} corresponding to both sampling methods were calculated, and a summary of the results is available both in Table~\ref{table:conditionRank} and in Figure~\ref{figure:Conditioning}, left.  We see that $\kappa_2$ for the spiral pattern have a reasonable behavior for low radial orders $n$, but grow large for higher degrees\footnote{One of the several possible explanations can lie in the intrinsic structure of the spiral points: the radial density of the sampling is almost constant when approaching the border of the disk. Thus, neither they correlate well with the increasingly oscillatory behavior of  the Zernike polynomials close to the boundary, nor they approximate appropriately the pluripotential theory equilibrium measure of the disk.}, being greater than $10^6$ for $n=15$. In comparison, observe that the maximum condition number for the OCS, corresponding to a radial order of 30, is less than 100. In order to stay within this frame with the spiral sampling, we must restrict ourselves to Zernike polynomials of radial order not greater than 6.

\begin{table}
\centering
 \begin{tabular}{|c|c|c|}
\hline  $\mathbf{n}$ $\mathbf{[N]}$
&  {\textbf{ $\kappa_2$ for spiral sampling}} &  {\textbf{ $\kappa_2$ for OCS}} \\   \hline
10 [66] & 	$5.8 \times 10^{4}$ 	  		& 3.2 	  \\ \hline
15 [136] & $2.4 \times 10^{6}$	 	& 5.7  	  \\ \hline
20 [231] & $2.7 \times 10^{8}$	  	& 11.3 	  \\ \hline
22 [276] &	 $7.0 \times 10^{15}$ 	& 15.2  	  \\ \hline
27 [406] & $1.9 \times 10^{16}$	  	& 32.8 	  \\ \hline
30 [496] & $3.0 \times 10^{17}$	 	& 53.3  	  \\ \hline
\end{tabular}
\caption{Dependence of the condition number  $\kappa_2$ on the maximal Zernike radial orders $n$ and on the total number of polynomials $N$, both for the spiral sampling and for the optimal concentric sampling (OCS).}
\label{table:conditionRank}
\end{table}

We know that the high condition number of the collocation matrix $A_n$ in \eqref{defA_n} has a direct impact on the accuracy of the solution of the interpolation problem. For illustration, a synthetic wavefront given by a vector of 496 randomly generated Zernike coefficients between -1 and 1 was  sampled according to both schemes. These samples were used as the interpolation data for recovery of the original Zernike coefficients by solving the system with the collocation matrix \eqref{defA_n}. This procedure was repeated 100 times and the mean and standard deviation of the RMS error distribution for the coefficients were computed, as depicted in Figure~\ref{figure:Conditioning}, right. 
The  RMS  of $10^{-5}$ can be considered as a threshold above which the recovered coefficients are unreliable (less than $5$ accurate decimal digits). Observe the correspondence between both graphs in Figure~\ref{figure:Conditioning}.

\begin{figure}
\centering 
\hspace{-5mm}\includegraphics[scale=0.5]{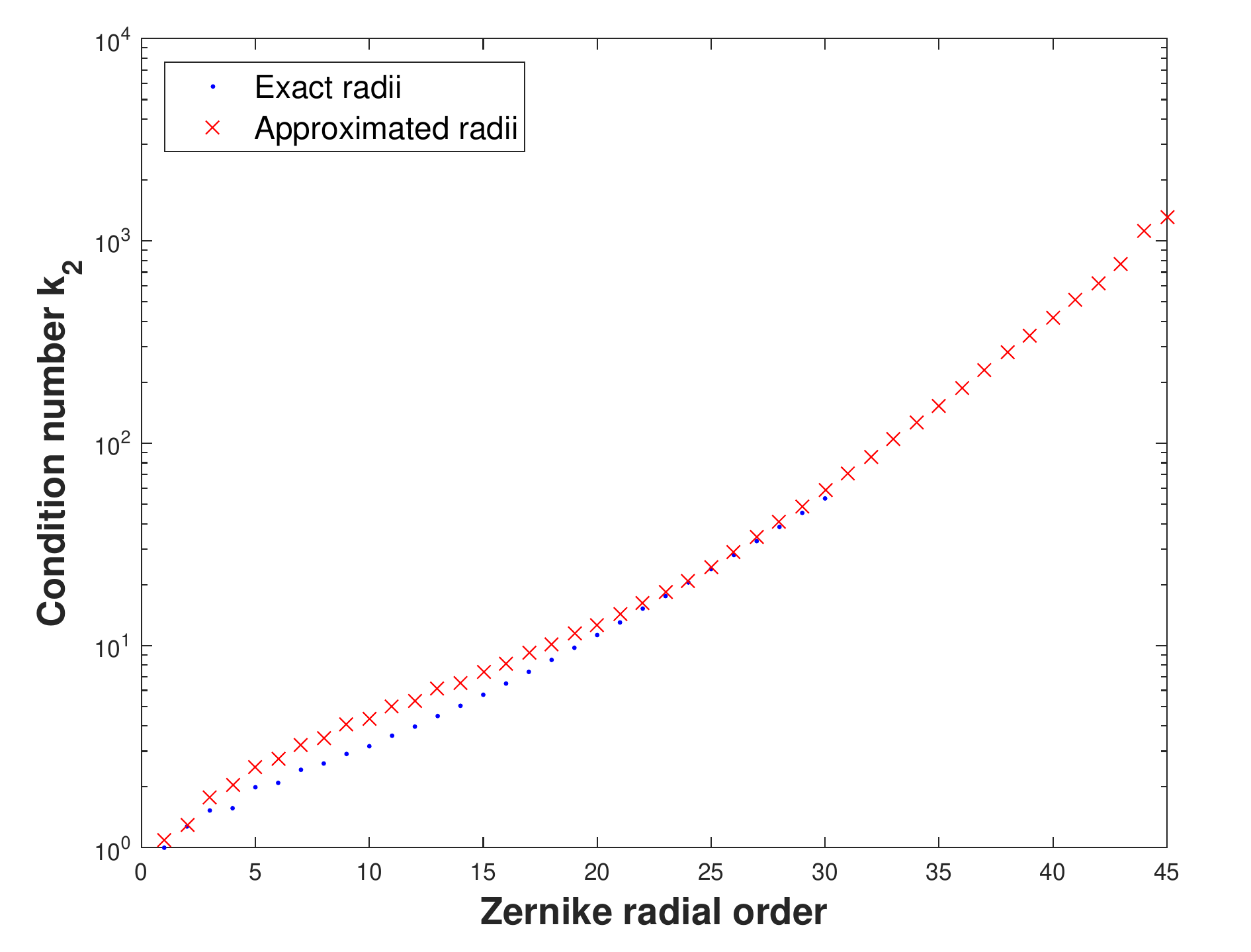}
 \caption {The effect of using the approximated values of the radii $r_i(n)$ given by \eqref{recCheb}, instead of the optimal ones computed numerically,  on the condition number of the matrix \eqref{defA_n}.}
\label{figure:factors}
\end{figure}

As it was explained, formula \eqref{recCheb} renders  an approximation to the optimal radii $r_i(n)$, found numerically. The impact of replacing the optimal values by this approximation on the numerical condition of the matrix \eqref{defA_n} can be appreciated in Figure~\ref{figure:factors}; we can see that the values of $r_i(n)$ are quite robust: the approximation hardly affects the condition number, at least for radial orders $\leq 30$ (we have computed numerically $r_i(n)$ for $n\leq 30$).

The nodes found in \cite{Carnicer} are prescribed by \eqref{nodes}, but the radii are given by the formula
$$
r_j=1-\left( \frac{2(j-1)}{n} \right)^a, \quad j=1, \dots, k(n),
$$
where the exponent $1<a<2$ is found experimentally and can depend on $n$; the value $a=1.46$ is said to give reasonable results for all degrees $n$. Numerical experiments show that these nodes  give moderate condition numbers $\kappa_\infty(A_n)$, although the nodes \eqref{recCheb}--\eqref{nodes} outperform them, using as the criterion the values of either $\kappa_2(A_n)$ or $\kappa_\infty(A_n)$.

Recall that the OCS interpolation nodes $\{P_i\}$ were built assuming that each ring has one node with argument $0$ (i.e., aligned with the positive $OX$ semi-axis), and then optimizing the radii $r_j$, see \eqref{recCheb}--\eqref{nodes}. A curious fact, already observed before (see e.g.~\cite{Cuyt2012,Carnicer}) is that the relative location of the equidistant nodes on each ring has a relatively low impact on the condition number of the collocation matrix. 
  
 The singular values of $A_n$, used for the computation of $\kappa_2$, can be computed as the square roots of the eigenvalues of the matrix $B=B_n=A_n^T A_n$, where $A_n^T$ is the matrix transpose of $A_n$. Observe that the elements $b_{ij}$ of $B=(b_{ij})$ are given by
$$
b_{ij} = \sum_{s=1}^N Z_{i-1}(P_s) Z_{j-1}(P_s) =\sum_{m=1}^k \sum_{\|P_s\|=r_m} Z_{i-1}(P_s) Z_{j-1}(P_s), 
$$
where $\|P\|$ is the distance of the node $P\in \R^2$ to the origin.  Observe that for $\ell=1, \dots, k(n)$,
$$
\sum_{\|P_s\|=r_\ell} Z_{i-1}(P_s) Z_{j-1}(P_s)=\gamma_{i-1}\gamma_{j-1} R_{i-1}(r_\ell) R_{j-1}(r_\ell)\Omega_{ij}(\ell),
$$
with
$$
\Omega_{ij}(\ell)= \sum_{s=1}^{n_\ell} f_i \left(|m_i | \theta^\ell_s(\alpha_{\ell})  \right) f_j \left(|m_i |\theta^\ell_s(\alpha_{\ell})  \right), \quad \theta^\ell_s(\alpha_{\ell})=\alpha_{\ell} +  \frac{2\pi(s-1)}{n_\ell}  ,
$$
where $n_\ell = 2 n+5-4\ell $, functions $f_{i}$ are either $\sin$ or $\cos$, and $0\leq \alpha_{\ell} < 2\pi/n_{\ell}$ is a parameter. The dependence of $m_i$ from $i$ is given by \eqref{formulaJ}. By standard trigonometric formulas,
$$
f_i \left(|m_i | \theta^\ell_s(\alpha_{\ell}) \right) f_j \left(|m_i | \theta^\ell_s(\alpha_{\ell})  \right) = \pm h \left((|m_i| - |m_j|)  \theta^\ell_s(\alpha_{\ell})  \right) \pm h \left( (|m_i| +|m_j|)  \theta^\ell_s(\alpha_{\ell})  \right) ,
$$
where $h=\frac{1}{2}\cos$ if $f_i=f_j$, and $h=\frac{1}{2}\sin$ otherwise. Consequently, 
$$
\Omega_{ij}(\ell)= \pm \sum_{s=1}^{n_\ell} h \left((|m_i| - |m_j|)  \theta^\ell_s(\alpha_{\ell})  \right) \pm \sum_{s=1}^{n_\ell} h \left((|m_i| + |m_j|)  \theta^\ell_s(\alpha_{\ell})  \right).
$$

In particular, let $\ell=1$;  values $||m_i| \pm |m_j||\leq 2n$ are not divisible by $n_1=2n+1$ except for $||m_i| \pm |m_j||=0$, and straightforward calculations show that 
$$
\sum_{s=1}^{n_1} h \left((|m_i| \pm |m_j|)  \theta^1_s(\alpha_{\ell})  \right) =\begin{cases}
\frac{n_1}{2}, & \text{if } m_i= m_j, \\
0, & \text{otherwise,}
\end{cases}
$$
and in any case, $\Omega_{ij}(1)$ is independent of $\alpha_{\ell}$. In other words,  the singular values of $A_n$ (and hence, the condition number $\kappa_2(A_n)$) do not depend on the exact location of the nodes $P_s$, equally spaced on the outermost ring.
 
\begin{figure}
\centering 
\includegraphics[width=0.75 \textwidth]{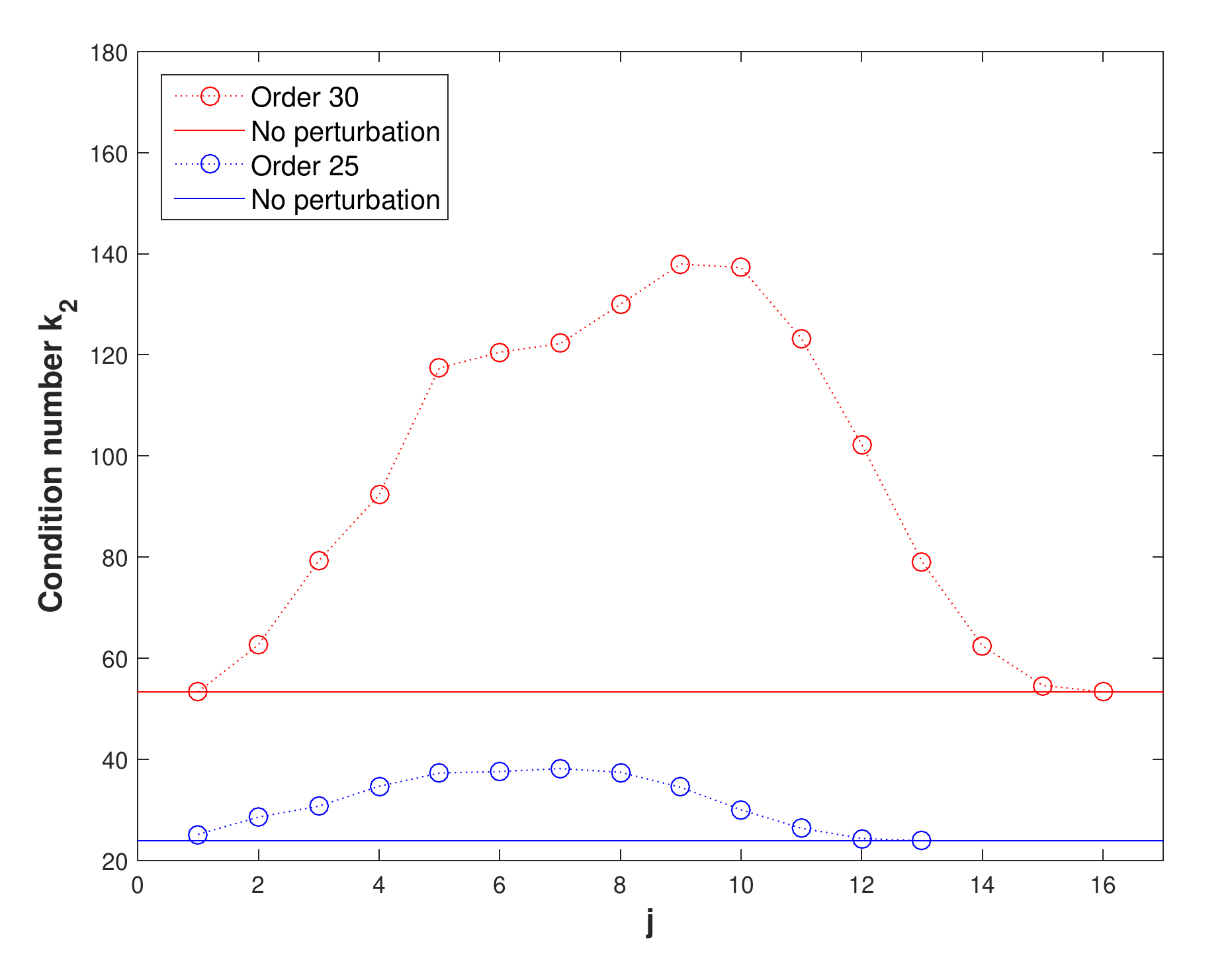}
 \caption {Maximal change in the condition number when a single circle of radius $r_j$ is rotated, $1\leq j \leq k(n)$, for radial orders $n=25$ and $n=30$.}
\label{figure:perturCircles}
\end{figure}

The considerations above are supported by some numerical experiments. First, for the maximal radial order of Zernike polynomials equal to $n=25$ (respectively, $n=30$), equally spaced nodes were placed on $k=13$ (respectively, $k=16$) concentric circles with radii given by~\eqref{recCheb}. Then the nodes on one of these circles were rotated preserving the equally spaced structure, and the maximal condition number $\kappa_2$ was computed. This procedure was repeated for each ring; the results are illustrated in Figure~\ref{figure:perturCircles}. 

\begin{figure}
\centering 
\includegraphics[width=0.75 \textwidth]{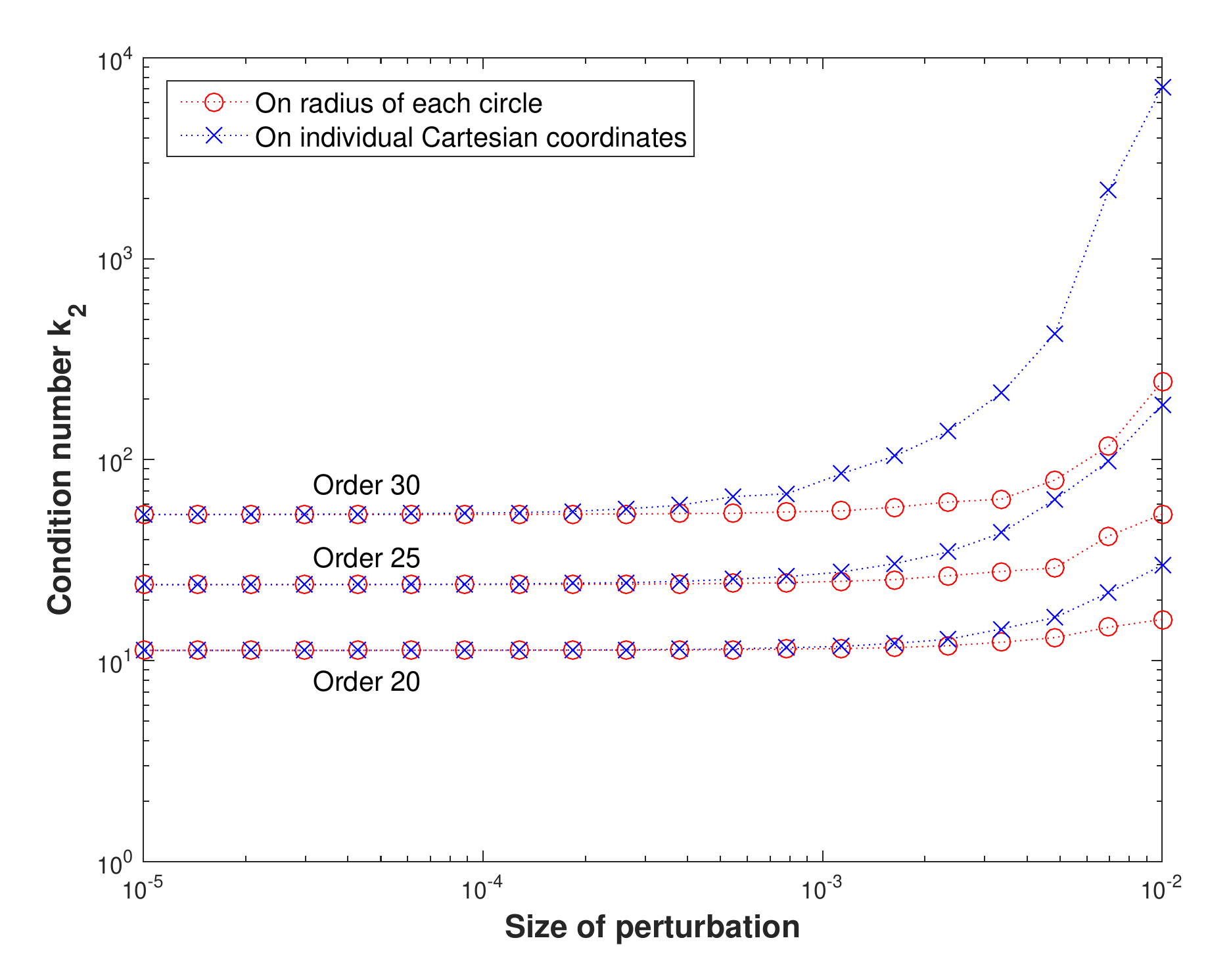}
 \caption {Change in the condition number when the radii of circles are perturbed ('o') and when individual nodes are perturbed ('x'), for radial orders 20, 25 and 30.}
\label{figure:perturTodo}
\end{figure}

Since in practice the precise placement of the interpolation nodes is difficult to guarantee, another criterion of usability of the OCS pattern is its stability to global perturbations of the circles and nodes. In Figure~\ref{figure:perturTodo} we illustrate the result of two independent experiments for radial orders $n=20, 25$ and $30$ of Zernike polynomials. The curves, marked with 'o', reflect the sensitivity of $\kappa_2$ to the variation of the optimal radii $r_j$, while the curves with 'x' show the change in $\kappa_2$ when individual nodes are randomly perturbed. It is worth pointing out that a study, similar in spirit, shows that weakly admissible meshes for the disk, unisolvent for polynomial interpolation, remain unisolvent and weakly admissible under small perturbations, see  \cite[Corollary 1]{Piazzon}.

\begin{figure}
\centering \begin{tabular}{lr}
\hspace{-10mm} \includegraphics[width=0.5 \textwidth]{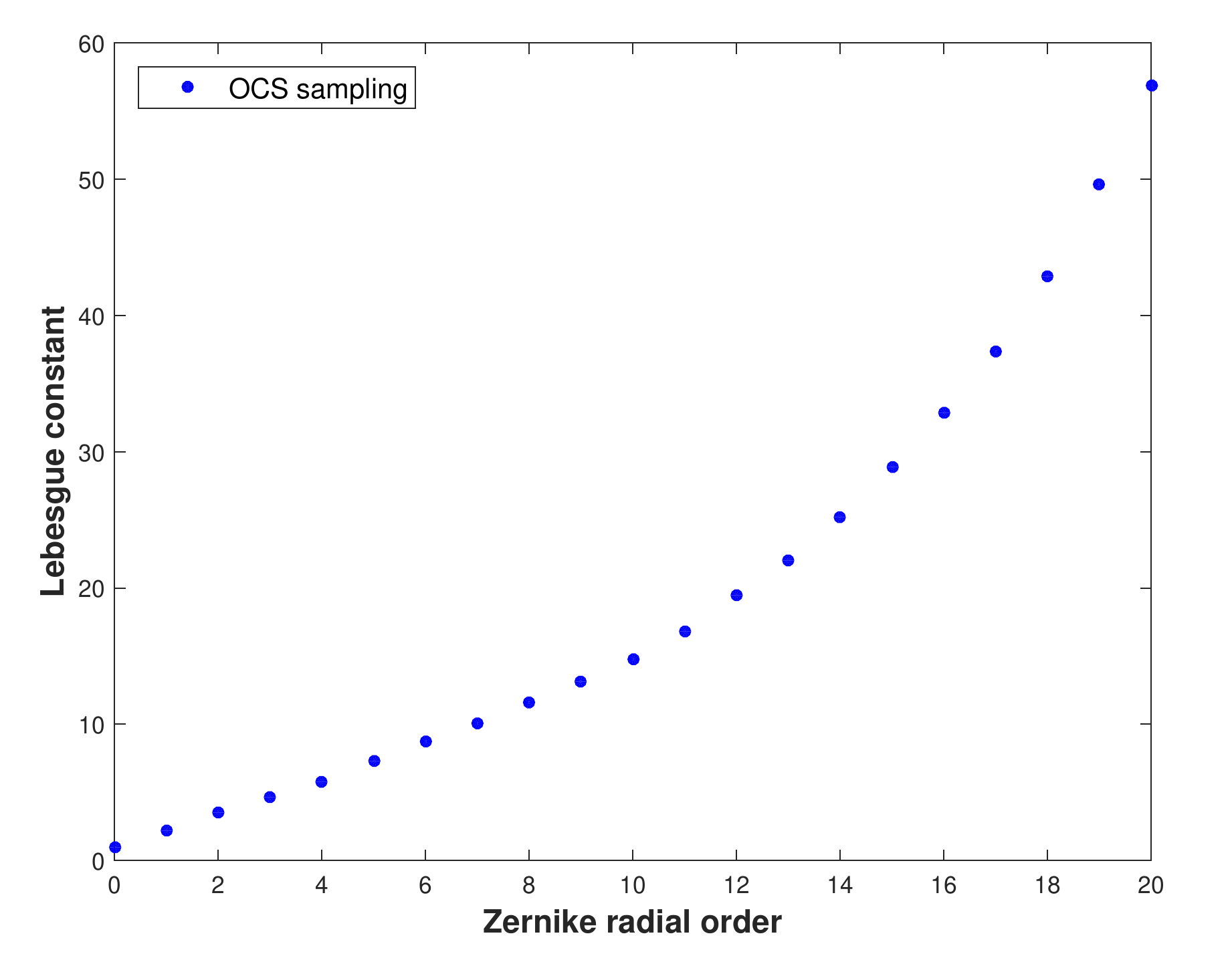} &
\hspace{-5mm} \includegraphics[width=0.5 \textwidth]{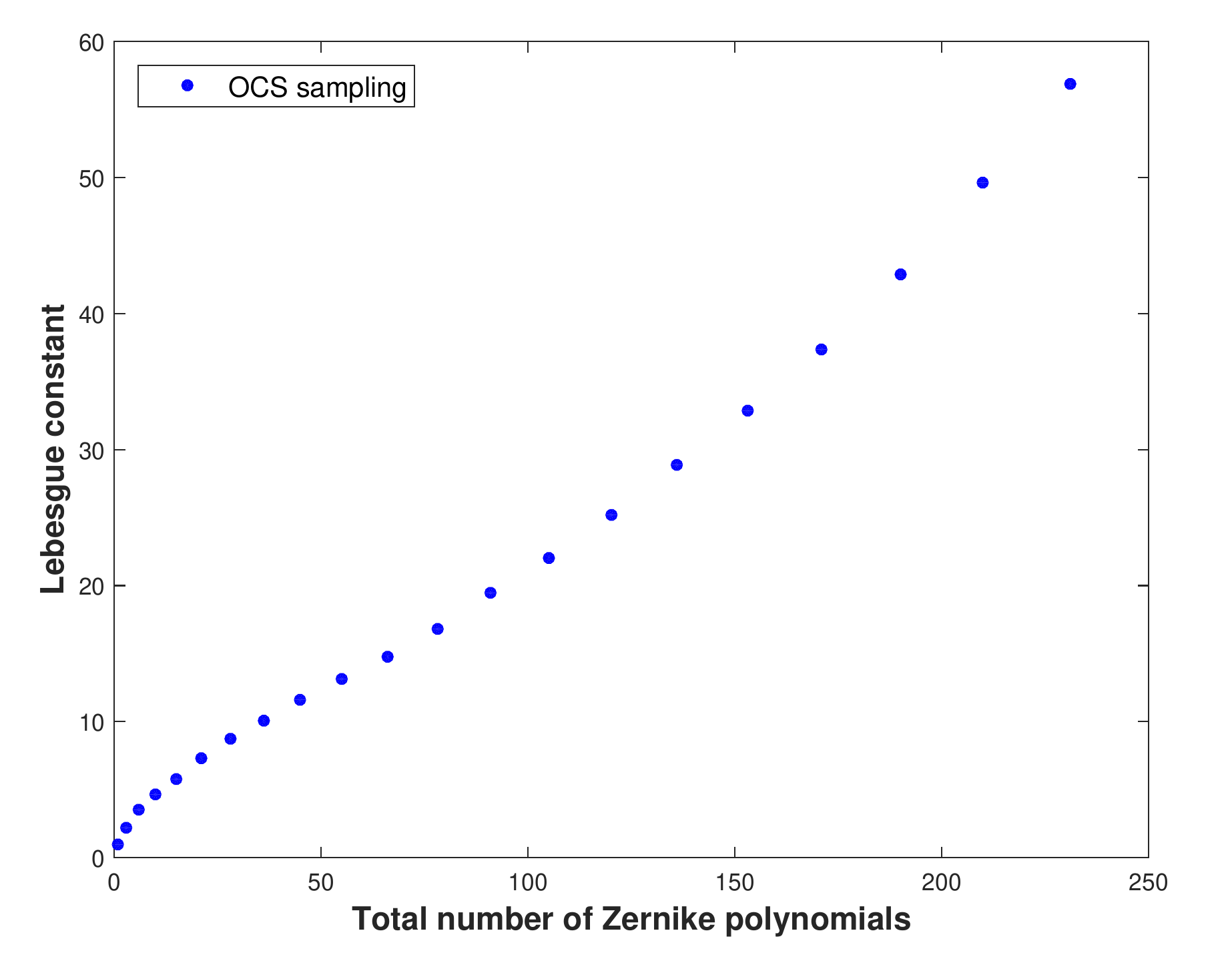}
\end{tabular}
 \caption {Lebesgue constants (vertical axis) as function of the order $n$ (left) and of the dimension of the interpolation space $N$ (right) for the case of the OCS.}
\label{figure:Lebesgue}
\end{figure}

We have also studied  the behavior of the Lebesgue constants corresponding to the OCS interpolation nodes $\{P_i\}$, using the formulas \eqref{lebesgue1}--\eqref{lebesgue2}. Obviously, straightforward numerical maximization of \eqref{lebesgue1} is a formidable task. A more efficient procedure was described recently in  \cite{Gunzburger}: it  uses an alternative way to assemble the Lagrange fundamental polynomials (following  \cite[Algorithm 4.1]{sauer}), along with replacing  maximization of $\ell(x,y)$ over the whole disk by its maximization on an admissible mesh. We have used both the ideas from  \cite{Gunzburger} and the direct evaluation of \eqref{lebesgue1} in a fine grid (definitely, a much slower procedure). The  results  are illustrated in Figure~\ref{figure:Lebesgue}. 

We see that in this case the dependence of $\Lambda_n$ from the dimension of the interpolation space $N$ is roughly linear, which is still far away from the  theoretically optimal $\mathcal O(\sqrt{n})$ from \cite{Sunder84}, but of the same order than the values obtained in \cite{Cuyt2012}. The Lebesgue numbers corresponding to the OCS interpolation nodes are also larger but comparable to those reported in 
\cite{Gunzburger} (see Table 3 and Figure 8 therein).

\begin{rmk}
Formula \eqref{recCheb} shows that the squares of the optimal radii are (asymptotically) uniformly distributed on $[0,1]$ with respect to the function
\begin{equation}
\label{defG}
G(x)=p\left(\sin \left(\frac{ \pi x}{2} \right)\right), \quad p(x)=\left(1.1565 \, x - 0.76535 \, x^2 + 0.60517 \, x^3\right)^2.
\end{equation}
$G$ is a  strictly increasing function on $[0,1]$, with $G(0)=0$, and $G(1)=0.992654<1$. 
As it follows from \cite{bloom92} (see also \cite[Section 8]{bloom12}), the interpolation nodes are asymptotically optimal (in the sense made more precise therein) if
$$
L(G):=\int_0^1 x^2 \log\left( G(x)\right) \, dx + 2 \int_0^1 \int_x^1 x \log\left( G(y)-G(x)\right)\, dy dx = - 2/3.
$$
Numerical integration shows that for $G$ given in \eqref{defG}, $L(G)= -0.681567$, which is slightly smaller. 

The two examples in \cite[Lemma 21 and 22]{bloom12} are
$$
G_1(x)=\sin^2 \left(\frac{ \pi x}{2} \right)\quad \text{and} \quad G_2(x)=1-(x^2-1)^2,
$$
for which $L(G_1)=-0.680609$ and $L(G_2)=-0.675676$, respectively. Notice that $G_1$ can be written in the form \eqref{defG} with  $p(x)= x^2$.
\end{rmk}

\begin{figure}
\centering 
\includegraphics[width=0.75 \textwidth]{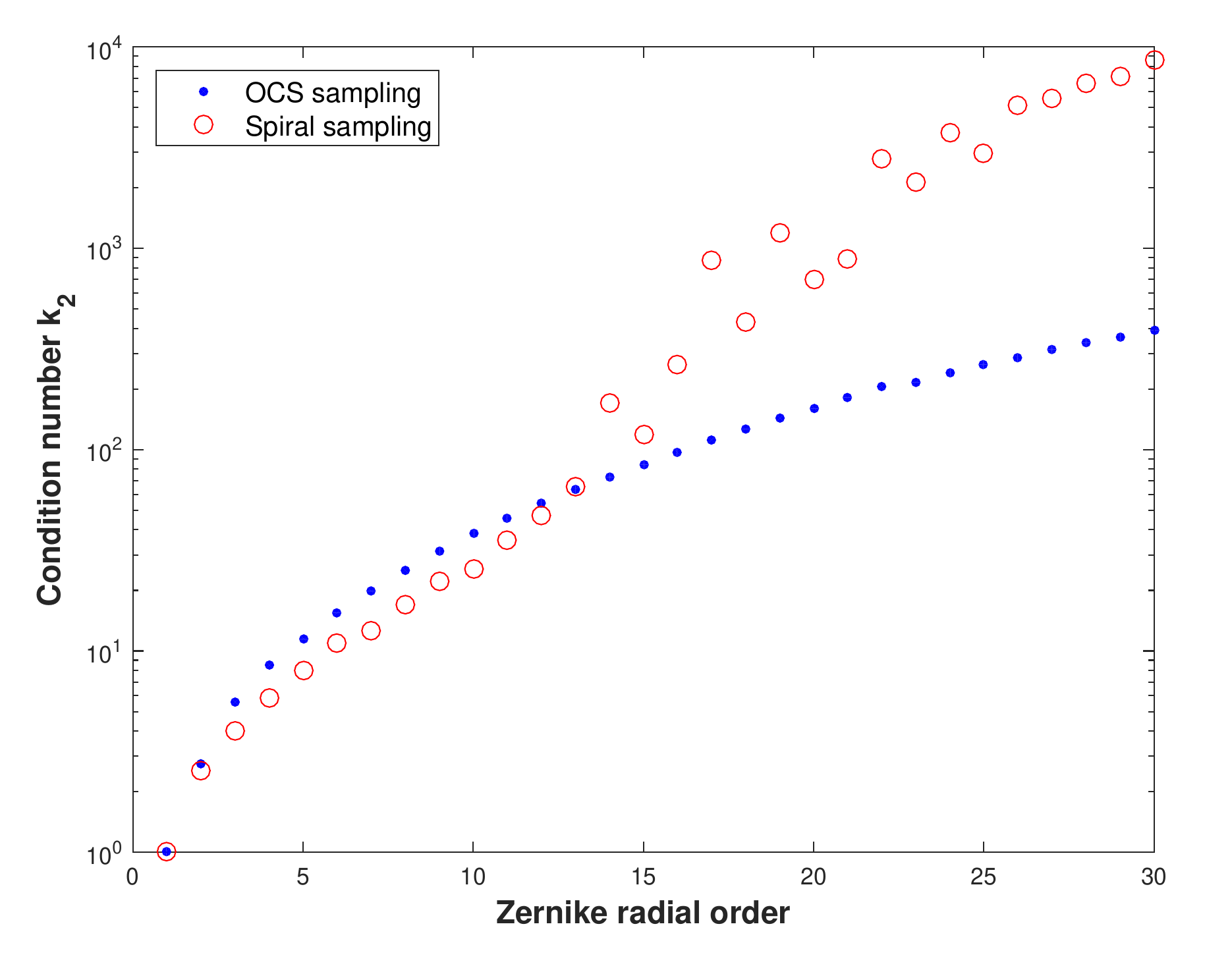} 
 \caption {Condition number of the collocation matrix, built from partial derivatives of the Zernike polynomials, as a function of the radial order, for both sampling patterns, namely, OCS and the spiral sampling.}
\label{figure:condDifferentials}
\end{figure}

As it was mentioned earlier, the wavefront sensors recover the actual wavefront from sampling its slopes at a finite number of points. Mathematically it boils down to solving systems of linear equations where the collocation matrix is built from the partial derivatives of the Zernike polynomials evaluated at the given nodes. This means, in particular, that each Zernike polynomial provides two rows of that matrix (corresponding to its partial derivatives), and that the first (constant) polynomial plays no role since its derivatives vanish. These considerations oblige to remove one sampling node for each pattern, in order to make the  size of the polynomial basis match the size of the sample. We chose to remove in each case the innermost node, although the results were very similar if any other sampling node was removed instead. The condition numbers of these collocation matrices of partial derivatives have been plotted in Figure~\ref{figure:condDifferentials} for both sampling patterns. 
In contrast to the results obtained for the matrix \eqref{defA_n}, now the condition numbers corresponding to the spiral and the OCS are similar, and in both cases   remain below $10^4$. We can observe also some  differences. For low orders (up to 13), the spiral sampling produces a relatively smaller condition number, although the difference is not significant. However, the condition number for the optimal concentric pattern grows slower with the radial order, so that for  radial orders above 13 we get better results than with  the spiral sampling. This difference can amount to about one order of magnitude for radial orders between 25 and 30.

\section{Conclusions}

The sampling nodes for polynomial interpolation, in particular using the Zernike basis, are crucial for an accurate reconstruction of an optical surface or wavefront from the sampled data. The optimal concentric sampling, defined by formulas \eqref{recCheb}--\eqref{nodes}, provides in many senses a quasi-optimal choice: this set is unisolvent (poised) for interpolation and exhibits a very moderate growth of both the condition numbers $\kappa_2$ and of the Lebesgue constants. This set is also fairly stable: as we have seen, small perturbations in the values of the radii $r_j$ or in the location of the individual nodes $P_j$ have no considerable influence on the condition number $\kappa_2$ of the collocation matrix. Since in general the precise placement of the interpolation nodes is difficult to guarantee, this fact is very relevant for the practical applicability of the interpolation scheme.

We have analyzed also the influence of rotation of the equally spaced nodes along each individual ring, showing that we can  choose the nodes on every circle independently from the others without affecting $\kappa_2$ considerably.
 
We can conclude that the optimal concentric sampling described in this paper not only renders a significant improvement in the accuracy of the recovery of the Zernike coefficients for low orders, but allows also  the possibility of using higher radial orders (of total degree 30 and even higher), which is not practical  with other sampling patterns.

\section*{Acknowledgements}

The first (DRL) and the fourth (AMF) authors were partially supported by MICINN of Spain and by the
European Regional Development Fund (ERDF) under grant
MTM2011-28952-C02-01, by Junta de Andaluc\'{\i}a (Excellence Grant P11-FQM-7276 and the research group
FQM-229), and by Campus de Excelencia Internacional del Mar (CEIMAR) of the
University of Almer\'{\i}a.  
This work was completed during a visit of AMF to the Department of Mathematics of the Vanderbilt University. He acknowledges the hospitality of the hosting department, as well as a partial support of the Spanish Ministry of Education, Culture and Sports through the travel grant PRX14/00037.

The second author (MASG) acknowledges the partial support of the Ministry of Economy and Competitiveness of Spain, Grant MTM2012-37894-C02-01.

The third author (MFM) especially acknowledges the valuable support provided by Centro Universitario de la Defensa en la Academia General del Aire de San Javier (Murcia, Spain).

We are also indebted to the referee of this paper for several bibliographical references and helpful comments.


\end{document}